\documentclass[11pt]{article}         %[charged.tex 29.06.08]
\usepackage{amsfonts,amssymb}
\makeatletter
\def\rdraft{\pagestyle{myheadings}
            \topmargin=30pt\headheight=10pt\headsep=3pt\footskip=20pt
            \textheight=10.5truein \textwidth=7.5truein
            \parindent=8pt \voffset=-1truein
            \ifcase \@ptsize \hoffset=-1.5truein \or \hoffset=-1.35truein
                    \or \hoffset=-1.15truein \fi}
\def\quality{\textheight=240mm \textwidth=160mm \topmargin=0Truein
             \ifcase \@ptsize \hoffset=-23mm
                     \or \hoffset=-20mm \or \hoffset=-15mm \fi}
\makeatother \quality %\rdraft

\def\beq#1#2{\begin{equation} \label{#1} #2 \end{equation}}
\def\bea#1{\begin{eqnarray*} #1 \end{eqnarray*}} \def\a{\!\!\!&\!\!\!\!&}
\def\n{\noindent}   \def\map{T}   
    
\def\IR{{\mathbb{R}}}  \def\IZ{{\mathbb{Z}}}  \def\cP{{\cal P}}
\def\toas#1{\stackrel{#1}{\longrightarrow}}
\def\ep{\varepsilon}  \def\phi{\varphi}   
  \def\mod1{\,({\rm mod\ } 1)\,}
\def\t{\tilde}   \def\dist{\varrho} \def\s{~\cdot}

\def\function#1{\left\{\!\!\!\begin{array}{ll} #1 \end{array} \right.}
\def\proof{\smallskip \noindent {\bf Proof. \ }}       %start of proof
\def\blanksquare{\,\,\,$\sqcup\!\!\!\!\sqcap$}         %blank  square
\def\qed{\hfill\blanksquare\linebreak\smallskip\par}   %end of proof

\def\thname{Theorem}     \def\lmname{Lemma}      \def\prname{Proposition}
\def\dfname{Definition}  \def\crname{Corollary}  \def\rmname{Remark}

\newtheorem{theorem}{\thname}[section]   %Numbering: Theorem--Other section
\newtheorem{lemma}{\lmname}[section]     %{lemma}[theorem]{Lemma}   section
\newtheorem{corollary}[lemma]{\crname}   %lemma

\newtheorem{dftn}{\dfname}[section]
\newtheorem{rmrk}[lemma]{\rmname}
\newenvironment{remark}{\begin{rmrk}\rm}{\end{rmrk}}     %lemma

             \catcode`@=11 \@addtoreset{equation}{section} \catcode`@=12

\newcommand\mlbscale{1pt} %to change: \renewcommand\mlbscale{1.3pt}
\newif\iffigs\figstrue %\newif\iffigs\figsfalse -- to fake figures
\def\bline(#1,#2)(#3,#4)(#5){\put(#1,#2){\line(#3,#4){#5}}}  %straight line

\def\bfig(#1,#2)#3#4{\begin{figure} \begin{center}
    \framebox{\setlength{\unitlength}{\mlbscale}
       \iffigs \begin{picture}(#1,#2) #3 \end{picture}
       \else \begin{picture}(60,10)(0,0)
                   \put(0,0){\framebox(60,10){Figure}} \end{picture} \fi}
    \end{center} \caption{#4} \end{figure}}
\def\Bfig(#1,#2)#3#4{\begin{figure} \begin{center}
    \setlength{\unitlength}{\mlbscale}
       \iffigs \begin{picture}(#1,#2) #3 \end{picture}
       \else \begin{picture}(60,10)(0,0)
                   \put(0,0){\framebox(60,10){Figure}} \end{picture} \fi
    \end{center} \caption{#4} \end{figure}}
\def\bpic(#1,#2)#3{\setlength{\unitlength}{\mlbscale}
    \begin{picture}(#1,#2) #3 \end{picture}}
    \def\V{V} 
 \def\*#1{#1^*} \def\Tr{\Delta}
\def\Free{{\rm Free}}  \def\Jam{{\rm Jam}}

\def\?#1{} % comments

\begin{document}%%%-------------------------------------------------%%%%
\title{Travelling with/against the flow. \\
       Deterministic diffusive driven systems.}
\author{Michael Blank\thanks{
        Russian Academy of Sci., Inst. for
        Information Transm. Problems,
        and Laboratoire Cassiopee UMR6202, CNRS, ~
        e-mail: blank@iitp.ru}
        \thanks{This research has been partially supported
                by Russian Foundation for Fundamental Research, CRDF
                and French Ministry of Education grants.}
       }
\date{10.10.08}
\maketitle
\begin{abstract} We introduce and study a deterministic lattice
model describing the motion of an infinite system of oppositely
charged particles under the action of a constant electric field.
As an application this model represents a traffic flow of cars
moving in opposite directions along a narrow road. Our main
results concern the Fundamental diagram of the system describing
the dependence of average particle velocities on their densities
and the Phase diagram describing the partition of the space of
particle configurations into regions having different qualitative
properties, which we identify with free, jammed and hysteresis
phases.
\end{abstract}

\n{\bf Keywords}: diffusive driven system, traffic flow, phase transition.

\section{Introduction}\label{s:intro}

It is well known that far from equilibrium statistical physics
systems may demonstrate a very complex behavior even on the level
of their stationary states. In general, such non-equilibrium
stationary (steady) states depend sensitively on the details of
the microscopic dynamics which leads to a very diverse
organization at the macroscopic level. An important class of such
systems is represented by {\em driven diffusive systems} (DDS),
describing the behavior of interacting diffusing particles driven
into selected spatial directions by an external force. Models of
this sort are widely used to describe vehicular and pedestrian
traffic \cite{MM}, water droplets in microemulsions with distinct
charges \cite{AN}, and numerous biological problems from molecular
motors to protein synthesis (see e.g. \cite{SZL}).

Basic ideas used in these models were introduced originally in
\cite{KLS} as a simple modification of the Ising lattice gas. From
that time tens of publications have been devoted to this subject.
From the point of view of mathematics these processes are close
relatives to asymmetric exclusion processes and topological Markov
chains.

The simplest model of this sort can be described as follows.
Consider a finite one-dimensional lattice with periodic boundary
conditions (i.e. a ring). Each site of the lattice is either vacant
($\cdot$) or occupied by a positive (+) or negative (-) particle.
Positive particles are driven to the right and negative particles
are driven to the left. The dynamics is random-sequential, i.e. at
each time step a pair of nearest neighbor sites is selected at
random and particles or vacancies located at these sites are
exchanged according to the following set of probabilities: %
$$\cP(+- ~\to~ +-)=\beta, ~~ \cP(-+ ~\to~ +-)=0, ~~ %
  \cP(+\s ~\to~ \s+)=\cP(\s- ~\to~ -\s)=1 .$$
Here the parameter $0<\beta\le1$ plays the role the time delay
required to exchange neighboring positive and negative particles
$+-$. Numerical studies and mean-field type approximations of this
model and its two-lane generalizations have been a subject of a
large number of publications (see e.g.
\cite{KLS,GLM,EKLM,KLMSW,KLMST,KLMT,GSZ,GSZ2,KSZ,ESMZ,AN,Ja} and
further references therein).

Our aim is to introduce a deterministic version of a DDS and to
study rigorously its limit dynamics (as time goes to infinity) on
the infinite integer lattice $\IZ$. In order to do this we first
pass from the random-sequential to synchronous movement of particles
(when all particles are trying to move simultaneously), and then
construct deterministic exchange rules for particles of different
signs. The simplest way to achieve the second goal is to replace the
random exchange with the introduction of the waiting time
$\tau=1/\beta$ before the exchange. (Indeed, the average waiting
time in the random exchange version is exactly $1/\beta$.) Further
on we assume that the parameter $\tau$ takes only positive integer
values.

To define the synchronous movement of particles one needs to resolve
the problem with the presence of triples $+\s-$ (positive and
negative particles separated by a single vacancy) additionally to
couples $+-$ considered in the asynchronous case. In the case of a
triple $+\s-$ both positive and negative particles are trying
simultaneously to exchange their positions with the same vacancy
(which cannot happen under the asynchronous updating). We assume
that in this case the waiting time before the particles exchange
their positions is equal to $\tau+1$.\footnote{
   One might argue that after the first time step both particles
   move to the vacant site and then use the standard waiting time
   $\tau$ to make further exchanges.} %
The situations $+-$ and $+\s-$ ~we call {\em short} and {\em long
interactions}, while the situations $++$ and $--$ when a particle
is blocked by another particle of the same sign we call {\em
simple interactions}. Using the notion of interactions the
dynamics can be described as follows. All particles in a
configuration are trying to move in the direction corresponding to
its sign by one position. If the particle does not interact with
others then it makes a move. If the short/long interaction takes
place then the particle waits for $\tau$ or $\tau+1$ time steps
(depending on the interaction type) before to start moving, while
in the case of the simple interaction the particle simply
preserves its position. As an example consider the dynamics of a
spatially periodic configuration with $\tau=2$ (only the main
period of length 17 containing 4 positive ({\Large\verb"+"}) and 2
negative ({\Large\verb"-"}) particles
and 11 vacancies ({\Large\verb"."}) is shown): %
\def\vb{\vskip-0.25cm}
\begin{center}{\Large
    \verb"+..+....-..++...-" \qquad {\small t=0}~~ \\ \vb
    \verb".+..+..-...+.+.-." \qquad {\small t=1}~~ \\ \vb
    \verb"..+..+-.....++.-." \qquad {\small t=2}~~ \\ \vb
    \verb"...+.+-.....++.-." \qquad {\small t=3}~~ \\ \vb
    \verb"....+-+.....+-.+." \qquad {\small t=4}~~ \\ \vb
    \verb"....+-.+....+-..+" \qquad {\small t=5}~~ \\ \vb
    \verb"+...-+..+...-+..." \qquad {\small t=6}~~ \\ \vb
    \verb".+.-..+..+.-..+.." \qquad {\small t=7}~~ \\ \vb
    \verb".+.-...+.+.-...+." \qquad {\small t=8}~~ \\ \vb
    \verb".+.-....++.-....+" \qquad {\small t=9}~~ \\ \vb
    \verb"+-.+....+-.+....." \qquad {\small t=10} \\ \vb
    \verb"+-..+...+-..+...." \qquad {\small t=11} \\ \vb
    \verb"-+...+..-+...+..." \qquad {\small t=12} \\ \vb
    \verb"..+...+-..+...+.-" \qquad {\small t=13} \\ \vb
    \verb"...+..+-...+..+.-" \qquad {\small t=14} \\ \vb
    \verb"....+.-+....+.+.-" \qquad {\small t=15} \\ \vb
}\end{center}

Here the first two positive particles and the first negative
particle move freely until the short interactions with the 1st
negative particle consecutively take place at time t=2 and t=4.
Then after 2 time steps the positive and negative particles
exchange their positions. The 3d positive particle at the 1st time
step makes a simple interaction with the 4th positive particle.
The 4th positive particle starting from t=1 is making a long
interaction with the 2nd negative particle during three time
steps, then exchange their positions, etc. Note that due to the
spatial periodicity after the 5th time step a positive particle
``moves'' from the 17th to the 1st position, while after the 12th
time step a negative particle ``moves'' from the 1st to the 17th
position.

In fact, this model is not Markov since the information about
present positions of particles alone does not allow to predict the
future motion of particles (due to the time delays during
particle interactions) and to make the model precise one needs to
introduce the notion of particles {\em states} which will be done
in Section~\ref{s:definitions}.

The analogy between the classical gas-liquid transition and
stationary states in interacting particle systems is well known and
discussed at length in the literature. The gaseous or laminar
phase corresponds to the situation when particles almost never
interact with each other and thus move with constant velocities.
Typically this happens at low particle densities. In the liquid
phase (which occurs normally at high densities) the particles are
located close to each other and mutual interactions organize them
into clusters (traffic jams). Additionally there might be an
intermediate or hysteresis phase when both types of particle
behavior coexist.

In our case the classical laminar phase cannot take place since
particles moving in opposite directions interact inevitably. To
overcome this difficulty we say that a particle is becoming {\em
free} starting from a certain time $t_0\ge0$ if for $\forall t\ge
t_0$ it does not take part in simple interactions. In other words,
after an initial transient period of length $t_0$ the particle
stops interacting with particles of the same sign. Then instead of
the laminar phase we consider a pair of {\em eventually free}
positive and negative phases defined as subsets $\Free_\pm^\infty$
of the set of all particle configurations as follows. A
configuration belongs to $\Free_+^\infty$ if each its positive
particle becomes free after a finite transient period, which might
depend on the particle. The set $\Free_-^\infty$ is defined in the
same way but for negative particles. Similarly instead of the
liquid phase we consider two {\em eventually jammed} (notation
$\Jam_\pm^\infty$) phases having the property that clusters of
particles of the corresponding sign are always present after some
transient period. The cases when the length of the transient
period is equal to zero we denote by $\Free_\pm$ and $\Jam_\pm$
respectively. Note that the sets of eventually free and eventually
jammed configurations of different signs have nonempty
intersections and that it is possible that an initial
configuration has no interacting particles but under dynamics they
will start interacting and form ``jams''.

%%%%%%%%%%%%%%%%%%%%%%%%%%%%%%%%%%%%%%%%%%%%%%%%%%%%%%%%%%%%%%%%%%
%% Fundamental diagram (+-)
\Bfig(290,120)
      {\footnotesize{
       %% picture (a)
       \put(0,0){\vector(1,0){150}} \put(0,0){\vector(0,1){120}}
       \thicklines
       \bezier{100}(0,60)(20,75)(40,80)  %phase 1
       \bezier{100}(50,70)(87,20)(125,0) %phase 3
       %\bezier{10}(46,80)(54,86)(58,75)  %phase 2
       %\bezier{10}(58,75)(64,55)(70,70)  %---"---
       %\put(48,78){$\cdot$} \put(50,73){$\cdot$} %phase 2
       %\put(52,81){$\cdot$} \put(54,70){$\cdot$} %---"---
       %\put(56,74){$\cdot$} \put(58,70){$\cdot$} %---"---
       %\put(60,60){$\cdot$} \put(62,65){$\cdot$} %---"---
       %\put(64,81){$\cdot$} \put(66,70){$\cdot$} %---"---
       \thinlines
       \bezier{10}(40,80)(45,81.5)(50,82) %continuation of phase 1
       %\bezier{10}(63,87.5)(67,80)(70,70) %continuation of phase 3
       \bezier{25}(40,95)(40,47)(40,0) %1st boundary (46,78)
       \bezier{25}(50,95)(50,47)(50,0) %2nd boundary
       \bezier{100}(0,105)(20,100)(40,90) %\t{V} (-)-particles
       \bezier{100}(50,20)(90,20)(125,20) %\t{V} (-)-particles
       \put(20,105){$\t{V}$} \put(110,25){$\t{V}$} %(-)-particles
       \put(35,-8){$\rho_c$} \put(48,-8){$\rho'_c$}
       \put(145,-8){$\rho$} \put(120,-8){1} \put(-10,20){$\frac1\tau$}
       \put(20,60){$V$} \put(60,60){$V$}
       \put(-10,115){$V$}   \put(-10,105){$1$} \put(-6,-6){$0$}
       %% picture (b)
       \put(170,0){\bpic(120,120){
         \bline(0,0)(1,0)(120)   \bline(0,0)(0,1)(120)
         \bline(0,120)(1,0)(120) \bline(120,0)(0,1)(120)
         \bline(40,0)(1,3)(40)   \bline(0,40)(3,1)(120)
         %\bline(36,0)(1,3)(21) \bline(34,0)(1,4)(17.2) %shade
         \bline(0,120)(1,-1)(120)
         \thicklines
         \bline(30,0)(1,4)(30)   \bline(0,30)(4,1)(120)
         \put(15,15){A} \put(15,75){B} \put(75,15){C}
         \put(15,-10){$\frac1{\tau+2}$}
         \put(35,-10){$\frac1{\tau+1}$}
         \put(-18,25){$\frac1{\tau+2}$}
         \put(-18,42){$\frac1{\tau+1}$}
         \put(50,125){$\frac\tau{\tau+2}$}
         \put(75,125){$\frac\tau{\tau+1}$}
         \put(122,57){$\frac\tau{\tau+2}$}
         \put(122,77){$\frac\tau{\tau+1}$}
         \put(115,-8){$\rho$}  \put(-8,115){$\t\rho$}
         \put(-6,-6){$0$}
         }}
      }}
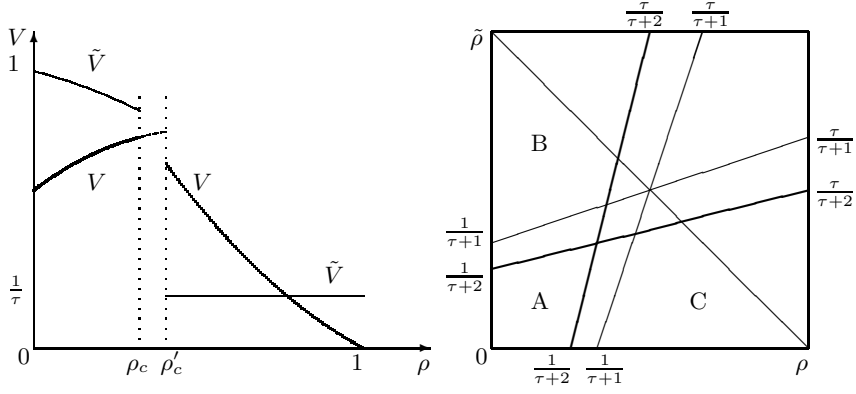
{(left) Fundamental diagram: dependence of the limit average
velocity $\V$ on the density $\rho$ of positive particles with the
fixed density $\t\rho<1/(\tau+2)$ of negative particles. The
average velocity $\t\V$ of negative particles is indicated by a
thin line. Values $\rho_c,\rho'_c$ indicate boundaries of the
hysteresis phase. %
\newline (right) Phase diagram:
$A\subset\Free_+^\infty\cap\Free_-^\infty$,
$B\subset\Free_+^\infty$, $C\subset\Free_-^\infty$. The region
$H:=\{0\le\rho+\t\rho\le1\}\setminus(A\cup B\cup C)$ between thick
and thin lines belongs to the hysteresis phase.
\label{f:fund-diag-pm}}
%%%%%%%%%%%%%%%%%%%%%%%%%%%%%%%%%%%%%%%%%%%%%%%%%%%%%%%%%%%%%%%%%

Our main results are formulated in the following two theorems.
We start with qualitative results justifying our discussion of
phase transitions and giving an inner characterization of
individual configurations belonging to each phase. To this end
in Section~\ref{s:BA} we introduce the notion of signed (positive
and negative) (proto)clusters\footnote{Roughly
   speaking a protocluster is a collection of particles
   of the same sign which will form a true cluster in future.
   A configuration $x\in\Jam_\pm^\infty$ may contain no infinite
   life-time clusters but in that case there are infinite
   life-time protoclusters.} %
of particles. By the life-time of a (proto)cluster one means the
duration of time until it ceases to exist. Exact definitions and
discussion see in Sections~\ref{s:definitions},\ref{s:BA}.

\begin{theorem}\label{t:phase-tr} Let for a configuration $x$ the
densities of positive and negative particles $(\rho(x),\t\rho(x))$
be well defined. If \\%
(a) $(\tau+2)\rho(x)<1+(\tau-1)\t\rho(x)$ then $x\in\Free_+^\infty$.\\
(b) $(\tau+2)\t\rho(x)<1+(\tau-1)\rho(x)$ then $x\in\Free_-^\infty$.\\
(c) $(\tau+1)\rho(x)>1+(\tau-1)\t\rho(x)$ then $x\in\Jam_+^\infty$.\\
(d) $(\tau+1)\t\rho(x)>1+(\tau-1)\rho(x)$ then
$x\in\Jam_-^\infty$.
\end{theorem}

To formulate quantitative results consider a partition of the
triangle $\Tr:=\{(\rho,\t\rho):~ 0\le\rho+\t\rho\le1,
                      ~ \rho,\t\rho\ge0\}$ %
(describing all possible pairs of densities) made by 4 straight
lines:
$$(\tau+1)\rho=1+(\tau-1)\t\rho,\quad
 (\tau+2)\rho=1+(\tau-1)\t\rho, $$
$$ (\tau+1)\t\rho=1+(\tau-1)\rho,\quad
 (\tau+2)\t\rho=1+(\tau-1)\rho $$ %
(see Fig~\ref{f:fund-diag-pm} (right)). Denote %
\bea{ \a A:=\{(\rho,\t\rho)\subset\Tr:~~
         (\tau+2)\rho<1+(\tau-1)\t\rho, ~~
         (\tau+2)\t\rho<1+(\tau-1)\rho \}, \\
   \a B:=\{(\rho,\t\rho)\subset\Tr:~~
         (\tau+2)\rho<1+(\tau-1)\t\rho, ~~
         (\tau+1)\t\rho>1+(\tau-1)\rho \}, \\
   \a C:=\{(\rho,\t\rho)\subset\Tr:~~
         (\tau+1)\rho>1+(\tau-1)\t\rho, ~~
         (\tau+2)\t\rho<1+(\tau-1)\rho \}, \\
   \a H:=\Tr\setminus(A\cup B\cup C). } %
Theorem~\ref{t:phase-tr} implies that
$A\subset\Free_+^\infty\cap\Free_-^\infty$,
$B\subset\Free_+^\infty\cap\Jam_-^\infty$,
$C\subset\Free_-^\infty\cap\Jam_+^\infty$.

\begin{theorem}\label{t:main-pm}
Let the densities of both positive and negative particles
$\rho,\t\rho$ in the initial configuration be well defined. Then
the corresponding average particle velocities $V,\t{V}$ are well
defined as well and depend only on the particle densities
according to the following relations:
$$\begin{array}{ll}
     V(\rho,\t\rho)=\frac{1+(\tau-1)(\rho-\t\rho)}
                        {1+(\tau-1)(\rho+\t\rho)},~~
  \t V(\rho,\t\rho)=\frac{1+(\tau-1)(\t\rho-\rho)}
                        {1+(\tau-1)(\rho+\t\rho)}
                   &\mbox{if } (\rho,\t\rho)\in A \\
     V(\rho,\t\rho)=\frac1\tau,~~
  \t V(\rho,\t\rho)=\frac1\tau~(\frac1{\t\rho}-1)
                   &\mbox{if } (\rho,\t\rho)\in B \\
     V(\rho,\t\rho)=\frac1\tau~(\frac1\rho-1),~~
  \t V(\rho,\t\rho)=\frac1\tau
                   &\mbox{if } (\rho,\t\rho)\in C
\end{array}$$
\end{theorem}

Note that Theorem~\ref{t:phase-tr} gives a qualitative description
of the asymptotic dynamics, while Theorem~\ref{t:main-pm}
describes it in quantitative terms. It is worth mention also that
without a systematic preliminary numerical modelling which allow
us to understand qualitatively the local structure of the particle
flow even the formulation of results proven in this paper would be
impossible.

Proofs of these results ideologically are based on the machinery
of the analysis of life-times of particle clusters developed in
\cite{Bl2,Bl3,Bl4} where deterministic interacting particle
systems with a single type of particles were studied. The
presence of particles moving in the opposite direction together
with two types of interactions (short and long) complicates
significantly the behavior of the system. In particular, new
clusters may be born in a free flow of particles, and the
asymptotic dynamics depends on two parameters (densities of
positive and negative particles). Therefore the technics has been
changed a lot and yet we are able to give only lower and upper
estimates of the life-times. Nevertheless these estimates allow
to find exact boundaries of all phases present in the model.

The paper is organized as follows. In Section~\ref{s:definitions}
we give the formal description of the model and main statistical
quantities under study: particle densities, average velocities,
etc. In the absence of clusters the dynamics of particles is
trivial (except from interactions between particles of opposite
signs). Therefore to study the model we need to analyze the
dynamics of clusters of particles and their ``life-times''. Exact
definitions of these objects and corresponding mathematical
results are discussed in Section~\ref{s:BA}. Duality relations
between positive and negative particles allow us to calculate in
Section~\ref{s:duality} average particle velocities under the
assumption that they are well defined. The latter is connected to
the proof of our main results formulated in the Introduction and
given in Section~\ref{s:proofs}. Section~\ref{s:hysteresis} is
dedicated to the analysis of the region of the Phase diagram where
different phases coexist. In Section~\ref{s:active} we apply our
results to study a model of an active tracer moving with or
against the particle flow, and Section~\ref{s:generalization} is
dedicated to the generalization of qualitative results for the
case of configurations for which particle densities are not well
defined.

\section{The model}\label{s:definitions}

The main disadvantage of the model discussed in the Introduction
is that the information about the present positions of particles
alone without the knowledge for how long currently occurring
interactions already take place does not allow to define the
future motion of particles. In order to overcome this difficulty
we introduce the notion of a state of a particle which takes into
account the complete information about the occurring interactions.
Let us give formal definitions.

By a {\em configuration} we mean a bi-infinite sequence
$x=\{x_i\}_{-\infty}^\infty$ with elements from the alphabet
$\{-\tau-1,-\tau,\dots,-1,0,1,\dots,\tau,\tau+1\}$. Positive
entries correspond to positive particles, negative entries to
negative ones, while zero entries correspond to vacancies. Non
zero entries will be referred as {\em states} of particles located
at corresponding sites. The states will be used to take into
account the delays during interactions of positive and negative
particles. Thus the largest value of the state $\tau+1$
corresponds to the long interaction.

The set of {\em admissible} configurations $X$ consists only of
configurations $x$ satisfying
the condition that for each $i\in\IZ$ %
\begin{itemize}
\item if $x_i>1$ then either $x_{i+1}=-x_i$ or $x_{i+1}=0$ and
$x_{i+2}=-x_i$, additionally $x_i=\tau+1$ implies $x_{i+1}=0$; %
\item if $x_i<-1$ then either $x_{i-1}=-x_i$ or $x_{i-1}=0$ and
$x_{i-2}=-x_i$, additionally $x_i=-\tau-1$ implies $x_{i-1}=0$.
\end{itemize}
These conditions imply restrictions only to positions of
interacting particles having states greater than one on modulus.
For example, the configuration $\dots110\t101\t100020\t200\dots$
is admissible (here $\t1$ and $\t2$ stand for $-1$ and $-2$),
while $\dots110\t101\t100020\t100\dots$ is not, since the last two
particles separated by a single vacancy $20\t1$ are supposed to be
mutually interacting (according to their positions) but their
states differ on modulus.

By $x^t:=\map^tx,~t\in\IZ_+\cup\{0\}$ we denote the state of the
configuration $x$ at time $t$, assuming that the initial state
$x^0$ is given by $x$. Here $\map:X\to X$ is the map describing
the dynamics which we define on the level of individual particles
in the configuration $x\in X$ in the following three steps:
\begin{enumerate}
\item First consider sites $i,i'$ with $|i'-i|\le2$ containing
      mutually interacting particles with $x_i^t>1$ and $x_{i'}^t=-x_i^t$
      and set $x_i^{t+1}:=x_i^t+1, ~~ x_{i'}^{t+1}:=x_{i'}^t-1$. Then
      if $i'-i=1$ (short interaction) and $x_i^{t+1}>\tau$ or
      if $i'-i=2$ (long interaction) and $x_i^{t+1}>\tau+1$ set
                 $x_i^{t+1}:=-1, ~~ x_{i'}^{t+1}:=1$.
\item Then consider the sites $i$ with $x_i^t=1$. \\
      (a) if $x_{i+1}^t=0$ and $x_{i+2}^t\ne-1$ then set
         $x_i^{t+1}:=0, ~ x_{i+1}^{t+1}:=1$;\\
      (b) if $x_{i+1}^t=0$ and $x_{i+2}^t=-1$ then set
         $x_i^{t+1}:=2, ~ x_{i+2}^{t+1}:=-2$;\\
      (c) otherwise if $x_{i+1}^t=-1$  then set
         $x_i^{t+1}:=2, ~ x_{i+1}^{t+1}:=-2$.

\item It remains to consider the sites $i$ with $x_i^t=-1$
      which were not taken into account during the step 2.
      If $x_{i-1}^t=1$ set $x_{i-1}^{t+1}:=-1, ~ x_{i}^{t+1}:=0$;
      otherwise do nothing.
\end{enumerate}

In words, if a particle is not interacting it simply moves by one
position in the direction corresponding to its sign (rule 2.a and
3). In case of the simple interaction (the particle is blocked by
another particle of the same sign) it does not move an does not
change its state. If the short/long interaction takes place the
particle preserves its position but its state changes by $\pm1$
depending on the particle sign (see rule 1) until it reaches on
modulus the value $\tau$ (in the case of the short interaction) or
$\tau+1$ (in the case of the long interaction). After that the
particles get the states $\pm1$ (preserving original signs) and
exchange their positions. The rules 2.b and 2.c take care about
the initial stage of interactions.

Using the notation $\t1=-1, \t2=-2, \t3=-3$ we can rewrite the
example of dynamics of a spatially periodic configuration (see the
formal definition below) with $\tau=2$ described in the
Introduction as follows: %
\bigskip %
\def\vb{$$\vskip-0.65cm$$}   {\large\bf $$
  1  0  0  1  0  0  0  0\t1  0  0  1  1  0  0  0\t1\qquad t=0~~\vb
  0  1  0  0  1  0  0\t1  0  0  0  1  0  1  0\t1  0\qquad t=1~~\vb
  0  0  1  0  0  1\t1  0  0  0  0  0  1  2  0\t2  0\qquad t=2~~\vb
  0  0  0  1  0  2\t2  0  0  0  0  0  1  3  0\t3  0\qquad t=3~~\vb
  0  0  0  0  1\t1  1  0  0  0  0  0  1\t1  0  1  0\qquad t=4~~\vb
  0  0  0  0  2\t2  0  1  0  0  0  0  2\t2  0  0  1\qquad t=5~~\vb
  1  0  0  0\t1  1  0  0  1  0  0  0\t1  1  0  0  0\qquad t=6~~\vb
  0  1  0\t1  0  0  1  0  0  1  0\t1  0  0\t1  0  0\qquad t=7~~\vb
  0  2  0\t2  0  0  0  1  0  2  0\t2  0  0  0  1  0\qquad t=8~~\vb
  0  3  0\t3  0  0  0  0  1  3  0\t3  0  0  0  0  1\qquad t=9~~\vb
  1\t1  0  1  0  0  0  0  1\t1  0  1  0  0  0  0  0\qquad t=10\vb
  2\t2  0  0  1  0  0  0  2\t2  0  0  1  0  0  0  0\qquad t=11\vb
\t1  1  0  0  0  1  0  0\t1  1  0  0  0  1  0  0  0\qquad t=12\vb
  0  0  1  0  0  0  1\t1  0  0  1  0  0  0  1  0\t1\qquad t=13\vb
  0  0  0  1  0  0  2\t2  0  0  0  1  0  0  2  0\t2\qquad t=14\vb
  0  0  0  0  1  0\t1  1  0  0  0  0  1  0  3  0\t3\qquad t=15\vb
$$}
%\bigskip

Invariance of the set of admissible configurations follows
immediately from the definition of the map $\map$. Observe also
that the restriction of the map $\map$ to the subset of admissible
configurations containing only entries from the alphabet $\{0,1\}$
coincides with the classical Nagel-Schreckenberg traffic flow
model (see e.g. \cite{NS,Bl3}).

For a configuration $x\in{X}$ by the {\em density} of positive
particles $\rho(x,I)$ in a finite lattice segment\footnote{We
  shall use also the notation $x[n,m]$ for lattice segments to
  specify the configuration $x$.} %
$I=[n,m]:=\{i\in\IZ:~~n\le{i}\le m\}$ we mean the number of
positive particles from the configuration $x$ located in $I$
divided by the total number of sites in $I$ (notation $|I|$), and
by $\t\rho(x,I)$ the corresponding value for the negative
particles. If for any sequence of {\em nested} finite lattice
segments $\{I_n\}$ with $|I_n|\toas{n\to\infty}\infty$ the limits
$$\rho(x):=\lim_{n\to\infty}\rho(x,I_n), \qquad
  \t\rho(x):=\lim_{n\to\infty}\t\rho(x,I_n)$$
are well defined and do not depend on $\{I_n\}$ we call $\rho(x)$
the {\em density} of positive particles and $\t\rho(x)$ the
density of negative particles in the configuration $x\in{X}$.
Otherwise one considers upper and lower particle densities
$\rho^\pm(x),~\t\rho^\pm(x)$.

\begin{lemma}\label{l:density-preservation}
Particle densities are conserved under dynamics, i.e.
$\rho^\pm(x^t),\t\rho^\pm(x^t)$ do not depend on $t<\infty$.
\end{lemma}
\proof For a given lattice segment $I\in\IZ$ the number of
particles from the configuration $x^t\in{X}$ which can leave it
during the next time step cannot exceed 2 and the number of
particles which can enter this segment also cannot exceed 2. A
close look shows that the total number of particles that can leave
or enter this segment cannot exceed 2, because if a particle
leaves the segment through one of its ends no other particle can
enter through the same end. Therefore %
$$|\rho(x^t,I) - \rho(x^{t+1},I)|\cdot|I|\le2 $$ %
which implies the claim in the case of positive particles. The
proof in the case of negative particles is exactly the same. \qed

%Observe that we use an important property which holds for the
%systems we consider here: particle conservation.

A configuration $x\in X$ is said to be {\em spatially periodic}
with {\em period}\ $L\in\IZ$ if $x_i=x_{i+L}~~\forall i\in\IZ$.
Such configurations represent an important class of admissible
configurations for which both positive and negative particle
densities are well defined. A direct check demonstrates that the
spatial periodicity and its period\footnote{The {\em minimal}\
   spatial period may decrease under dynamics.} %
are preserved under dynamics. Therefore the dynamics of spatially
periodic configurations of period $L$ is equivalent to the
dynamics on the finite segment of the integer lattice of length
$L$ with periodic boundary conditions (i.e. to the only situation
suitable for numerical modelling). The total number of admissible
configurations on this lattice segment cannot exceed
$L^{2\tau+3}$, which implies eventual periodicity in
time\footnote{A trajectory becomes periodic in time after some
   finite transient period.} %
of the dynamics. Thus a trajectory $\{\map^tx\}_{t\ge0}$ of the
original infinite system with a spatially periodic initial
configuration $x\in X$ is eventually periodic in time as well.

\bigskip

We say that particles of the same sign located at sites $i'<i$ are
{\em consecutive} if all sites between them are either vacant or
occupied by particles of the opposite sign. Consider a pair of
consecutive positive particles. The following result shows that the
distance between these particles (calculated as $i-i'-1$) can shrink
to zero with time but cannot be enlarged much. (Compare with the
incompressibility property of a conventional fluid!)

\begin{lemma}\label{l:distance-pm} Let $\hat{X}$ be a
collection of admissible configurations $x$ having a pair of
consecutive positive particles at sites $i'<i$. Denote by $i'(x^t),
i(x^t)$ positions of these particles at time $t$ in the
configuration $x^t$ and by $D(x^t):=i(x^t)-i'(x^t)-1$ the
distance between them. Then %
\beq{e:dist-particles}{
 0 = \inf_{x\in\hat{X}} \liminf_{t\to\infty}D(x^t)
   < \sup_{x\in\hat{X}} \limsup_{t\to\infty}D(x^t)
   < 2\tau (i-i') .}%
\end{lemma}
\proof The lower estimate follows from the observation that,
independently on the initial distance between the particles, there
might be a large enough ``jam'' of positive particles ahead of them
which will stop the leading particle and allow the rear one to catch
up with it.

To prove the upper estimate observe that without interactions with
negative particles the distance $D(x^t)$ cannot grow in time. Note
also that a finite number of negative particles located initially in
the segment $[i',i]$ can give only a constant contribution to the
variation of the distance between the particles. Therefore it is
enough to consider only the case when initially there are only
vacancies in the segment $[i'+1,i-1]$. Denote by $t_1$ the first
moment of time when the leading particle meets a negative one, by
$t_1'$ the duration of the interaction between them (which might
take values $\tau$ or $\tau+1$), by $t_2$ the time between the end
of this interaction and the moment when the rear particle meets with
the same negative particle, and by $t_2'$
-- the duration of the latter interaction. Then we have %
\bea{
  i'(x^{t_1+t_1'})\a=\max\{i'(x^{t_1}) + t_1', i(x^{t_1})-1\}, \cr
    i(x^{t_1+t_1'})\a=i(x^{t_1}) + 1 + (t_1'-\tau) ;\cr
  i'(x^{t_1+t_1'+t_2})\a=i'(x^{t_1+t_1'}) + t_2 ,\cr
    i(x^{t_1+t_1'+t_2})\a\le i(x^{t_1}) + 1 + (t_1'-\tau) + t_2 ;\cr
  i'(x^{t_1+t_1'+t_2+t_2'})\a=i'(x^{t_1+t_1'+t_2})
                              + 1 + (t_2'-\tau), \cr
  i(x^{t_1+t_1'+t_2+t_2'})\a\le i(x^{t_1})
             + 1 + (t_1'-\tau) + t_2 + t_2'.}%
Therefore %
\bea{
   \a D(x^{t_1+t_1'}) = i(x^{t_1+t_1'}) - i'(x^{t_1+t_1'}) - 1 \cr
   \a\qquad = \max\{i(x^{t_1})-i'(x^{t_1}) - \tau + 1 -1, 2\}
           \le \max\{D(x^{t_1}),2\} ,\cr
\a D(x^{t_1+t_1'+t_2})
   = i(x^{t_1+t_1'+t_2}) - i'(x^{t_1+t_1'+t_2}) - 1 \cr
   \a\qquad \le i(x^{t_1}) + (t_1'-\tau) + t_2
     - \max\{i'(x^{t_1}) + t_1', i(x^{t_1})-1\} - t_2 \cr
   \a\qquad \le \max\{i(x^{t_1})-i'(x^{t_1}) -\tau-2, 2\}
           \le \max\{D(x^{t_1}),2\} ,\cr
\a D(x^{t_1+t_1'+t_2+t_2'})
   = i(x^{t_1+t_1'+t_2+t_2'}) - i'(x^{t_1+t_1'+t_2+t_2'}) - 1 \cr
   \a\qquad \le i(x^{t_1}) + (t_1'-\tau) + t_2 + t_2'
    - \max\{i'(x^{t_1}) + t_1', i(x^{t_1})-1\}
    - t_2 - 1 - (t_2'-\tau) \cr
   \a\qquad = i(x^{t_1}) + t_1' -1
    - \max\{i'(x^{t_1}) + t_1', i(x^{t_1})-2\} \cr
   \a\qquad \le \max\{i(x^{t_1})-i'(x^{t_1})-1, t_1'\}
           \le \max\{D(x^{t_1}), \tau+1\} .}%
This implies the result since the contribution from a single
negative particle located initially in the segment $[i',i]$ cannot
exceed $\tau+1$ and the number of such particles cannot be larger
than $i-i'-1$. \qed

\begin{remark} Result of Lemma~\ref{l:distance-pm} is not obvious
and is not what one expects here. Indeed, due to the presence of two
types of interactions: short and long, one would expect that if the
leading particle in the pair takes part only in short interactions
while another one takes part only in long interactions, then the
distance between them grows linearly in time.
Lemma~\ref{l:distance-pm} demonstrates that this is not the case. An
important property that we use here is that the difference between
the duration of long and short interactions is exactly equal to one.
If this would not be the case the distance between the particles
indeed might grow with time.
\end{remark}

Let us introduce now the notion of the {\em average velocity} of a
particle. For a configuration $x$ denote by $i^t$ the location of
a certain particle at time $t\ge0$ located originally at site
$i=i^0$. By a {\em finite time velocity} of this particle we mean
$V(x,i,t) := \frac1t~|i^t-i^0|$. If the limit
$V(x,i):=\lim_{t\to\infty} V(x,i,t)$ is well defined we call it
the {\em average velocity} of the particle. To distinguish between
the positive and negative particles in the latter case we use the
notation $\t{V}(x,i)$.

Our aim now is to show that the average velocity does not depend
on the choice of a particle.

\begin{lemma}\label{l:velocity-inv-pm}
Let $x$ be an admissible configuration and assume that for a
positive particle originally located at site $i$ the average
velocity $V(x,i)$ is well defined. Then for any positive particle in
the configuration $x$ the average velocity is well defined and
coincides with $V(x,i)$.
\end{lemma}
\proof Denote by $\hat{i}$ the location at time $t=0$ of the next
positive particle located to the right from $i$. As usual we
denote by $i^t, \hat{i}^t,~t\ge0$ positions of these particles at
time $t\ge0$. Then we have %
$$ V(x,\hat{i},t) = \frac1t~(\hat{i}^t - \hat{i}^0)
                = \frac1t~(i^t - i^0)
                + \frac1t~(\hat{i}^t - i^t)
                + \frac1t~(i^0 - \hat{i}^0) .$$
Applying Lemma~\ref{l:distance-pm}, according to which
$1\le\hat{i}^t - i^t\le2\tau|i - \hat{i}|$, we get
$$ |V(x,\hat{i},t) - V(x,i,t)|
   \le \frac1t~(\hat{i}^t - i^t) + \frac1t~|i^0 - \hat{i}^0|
   \le \frac{2\tau+1}t~|i - \hat{i}| \toas{t\to\infty}0 .$$
Thus $V(x,i)=V(x,\hat{i})$. Using the same argument one extends
this result to neighboring positive particles, and repeating it to
all positive particles in the configuration. \qed

\begin{corollary}\label{c:velocity-inv-pm} The invariance
of velocities of negative particles can be proven along the same
argument. Using these results we may drop the dependence on the
index $i$ in the definition of the average velocity.
\end{corollary}

\section{Dynamics of (proto)clusters and their basins of
attraction (BA)} \label{s:BA} %
For a collection of consecutive particles of the same sign ordered
with respect to their sign (i.e. according to their positions in the
case of the positive sign, and in the opposite way in the case of
the negative sign) we call {\em rear/leading} the first/last
particle with respect to this order. By a {\em cluster} (of
particles) we mean a segment $J=x[m',m], ~m'<m$ of the configuration
$x$ consisting only of particles of the same sign (called the {\em
body} of the cluster) and in which the leading particle is
interacting with a particle of the opposite sign (called the {\em
root} of the cluster) and the site next to the rear particle is
vacant or is occupied by a particle of the same sign as the root.
The sign of a cluster is identified with the sign of particles in
its body. Depending on the type of the interaction between the
leading particle and the root we say that the cluster has a {\em
short} or {\em long} root. The right hand side of the diagram below
shows positive clusters with short and long roots, while the left
hand side demonstrates how these clusters may be created:
$$\dots\underbrace{+\cdots+}_{{\rm body}}\cdots
     \underbrace{-}_{{\rm root}}\dots \quad \Longrightarrow \quad
  \function{\dots\underbrace{++}_{{\rm body}}
                 \underbrace{-}_{{\rm short~root}}\dots \\
            \dots\underbrace{++}_{{\rm body}}\cdot
                 \underbrace{-}_{{\rm long~root}}\dots}$$

In distinction to other deterministic transport flow models (see
e.g. \cite{Bl1,Bl2,Bl3,Bl4}) in the model under consideration new
clusters might be created even in the absence of clusters in the
initial configuration.

In order to control events of this sort we introduce the notion of
a {\em protocluster}\footnote{In the analysis of the free phase we
   will need to prove that the conditions of
   Lemma~\ref{l:life-time-upper-pm} imply that a flow of
   noninteracting at time $t\ge0$ positive particles remain free
   under dynamics.}
as a pair of consecutive particles of the same sign and a particle
of the opposite sign satisfying the property that these three
particles will form a true cluster and the moment when the leading
of these two particles will start interacting with the rear one is
the moment of the creation of a new cluster. In other words, they do
not join an existing cluster. The minimal lattice segment containing
this pair of particles is called the {\em body} of a protocluster
and the particle of the opposite sign with which they will form the
cluster -- its {\em root}. The sign of a protocluster is again
identified with the sign of particles in its body. Example of a
positive protocluster and clusters which it may form are shown in
the diagram above.
%$$\dots\underbrace{+\cdots+}_{{\rm body}}\cdots
%     \underbrace{-}_{{\rm root}}\dots \quad \Longrightarrow \quad
%  \dots\underbrace{++}_{{\rm body}}\underbrace{-}_{{\rm
%  root}}\dots$$
Note that between the body and the root of a protocluster there
might be some other particles of the same sign as the body, which
even might form a cluster with the root but will stop interacting
with the root particle before the leading particle of the
protocluster will reach them.

\begin{remark}\label{r:protocluster} A simple calculation shows
that particles of the same sign as the root might be present between
the body of the protocluster and its root only if they will make
short interactions with the particles from the body while with the
root they will make long interactions. Indeed, only in this case a
new cluster can be created.
\end{remark}

The boundaries of the body of a (proto)cluster might change with
time when the particles are moving (in the case of a protocluster)
or the leading particle leaves it while some particles join the
cluster from the other side. At time $t\ge0$ by $J^t$ we denote the
segment of the configuration corresponding to the body of the
(proto)cluster originally located at $J:=J^0$. If $t=0$ we drop the
time index to simplify the notation. The duration of time during
which the cluster exists (i.e. the number of particles in its body
is not smaller than two and the leading one is interacting with the
same root particle) we call the {\em life-time} of the cluster $J^t$
and denote by $LT(J)$. By the {\em life-time} of a protocluster we
mean the duration of time until the cluster formed by particles from
the body of the protocluster ceases to exist. If during its
life-time the body of a (proto)cluster consists of only two
particles (being present through its entire life-time) we call it
{\em trivial} and {\em nontrivial} in the opposite case. Clearly the
life-time of a trivial (proto)cluster is finite and depends only on
the distance to the root. Therefore we shall be interested only in
the case of nontrivial (proto)clusters. Observe that the situation
described in Remark~\ref{r:protocluster} corresponds to the trivial
protocluster.

`Attracting' preceding particles, a (proto)cluster $J$ plays a role
similar to an attractor in dynamical systems theory. Therefore it is
reasonable to study it in a similar way and to introduce the notion
of its {\em basin of attraction} (notation BA$(J)$), by which we
mean the minimal segment of the configuration $x$ containing all
sites from where particles may eventually join the (proto)cluster
during its life-time.

Due to the duality between positive and negative particles it is
enough to consider only positive (proto)clusters.  Let
$J=x[m',m],~m'<m$ be the body of a positive (proto)cluster having
its root at site $m''>m$. To take into account the dynamics we use
the notation $J^t=x^t[m'(t),m(t)]$ and $m''(t)>m(t)$ for the
corresponding objects at time $0\le t\le$LT$(J)$. We introduce two
functionals depending on an integer parameter $k<m'$: %
\beq{e:W-pm}{
   W(x[k,m]) := q(x_m)+\sum_{i=k}^m w(x_i), \quad
  W'(x[k,m]) := q(x_m)+\sum_{i=k}^m w'(x_i) ,}%
where
$$ w(z) := \function{\tau+1 &\mbox{if } z>0 \\
                     -\tau  &\mbox{if } z<0 \\
                     -1     &\mbox{otherwise } }, \qquad
  w'(z) := \function{\tau   &\mbox{if } z>0 \\
                     -\tau  &\mbox{if } z<0 \\
                     -1     &\mbox{otherwise } } ,$$
and $q(z)$ is defined as the amount of time until the leading
particle leaves the cluster, or $q(z):=\tau$ in the case of a true
protocluster. Thus $q(z)=\tau-z+1$ for the cluster with the short
root, and $q(z)=\tau-z+2$ for the cluster with the long root.

Denote by $j(t)$ the largest integer non-exceeding $m'(t)$ for
which $W(x^t[j(t),m(t)])=0$ or $W(x^t[j(t)-1,m(t)])<0$, and by
$j'(t)$ the largest integer non-exceeding $m'(t)$ for which
$W(x^t[j'(t),m(t)])=0$ or $W(x^t[j'(t)-1,m(t)])<0$. Note that both
$j(t)$ and $j'(t)$ might take an infinite value.

The following result allows to control the left boundary of a BA.

\begin{theorem}\label{t:BA-pm}
$x[j'(t),m(t)]\subseteq$BA$(J^t)\subseteq x[j(t),m(t)]$ and
$j(t+1)\ge j(t)+1,~ j'(t+1)=j'(t)+1$.
\end{theorem}

The proof of this Theorem will be given through a series of Lemmas,
but before to do this consider an example demonstrating that the
boundaries $j,j'$ might differ a lot and the quantity $j'-j$ even
might take an infinite value. Consider a spatially periodic
configuration $x$ with the main period containing three particles
and three vacancies: $+\cdot\cdot+\cdot~-$. These three particles
represent a protocluster for which $j(0)=-\infty, ~j'(0)=i-2$ where
$i$ stands for the position of the 1st positive particle.

\begin{lemma}\label{l:life-time-pm} The life-time of a cluster
of positive particles is equal to the number of positive particles
in its BA multiplied by $\tau$ minus the state of the leading
particle plus one in the case of the short root or plus two in the
case of the long root.
\end{lemma}
\proof The term ``the number of positive particles in its BA
multiplied by $\tau$'' describes the amount of time before all
particles in the BA will leave the cluster, assuming that the
initial state of the leading particle is one. Indeed, during the
life-time of a cluster the current leading particle is making
continuously the short interaction with the root particle. After
this the role of the leading particle is going to the next particle
in the original BA. The term ``the state of the leading particle
plus one ...'' takes into account the actual state of the leading
particle and the type of interactions with the root. \qed

\begin{lemma}\label{l:life-time-upper-pm}
BA(J)$\subseteq x[j,m]$ and the life-time of the cluster $J$ can
be estimated from above as the number of positive particles in
$x[j,m]$ multiplied by $(\tau+1)$.
\end{lemma}
\proof Consider the segment $I:=x[j,m]$, set $n:=m-j+1$, and denote
by $N$ the number of positive particles in this segment (which
cannot be smaller than 2 by the definition of the cluster) and by
$K$ the number of vacancies in it. Then the number of negative
particles in this segment is equal to $n+1-N-K$. Assume contrary to
our claim that there exists a particle $\xi$ located initially to
the left from the site $j$ which joins the cluster $J:=x[m',m]$
during its life-time. Since both vacancies and negative particles
are moving to the left (opposite to the movement of positive
particles) the particle $\xi$ must interact with all these vacancies
and negative particles before to join the cluster $J$. The amount of
time necessary in order to do this is at least
$t':=K+(n+1-N-K)\tau$. By the definition of the functional
$W(x[j,m])$ and the parameter $j$ we have
$$ N(\tau+1) - K - (n+1-N-K)\tau \le 0 .$$
Using this inequality we estimate the time $t'$ from above as
$N(\tau+1)$, which is larger than the upper estimate $N\tau+1$ of
the amount time during which all $N$ positive particles will join
the cluster $J$ and then leave it. Therefore the particle $\xi$
cannot join the cluster $J$ in time. We came to the contradiction.
\qed

\begin{lemma}\label{l:life-time-proto-pm}
Let $J$ be a nontrivial protocluster, then %
BA$(J^t)\subseteq x[j(t),m(t)]$.
\end{lemma}
\proof Denote by $\xi(t)$ the position at time $t\ge0$ of a
particle being originally at time $t=0$ to the left from the site
$j(0)$. Our aim is to show that $\xi(t)<j(t)$ for any moment of
time $t\ge0$ until the cluster is formed. Whence it is formed we
can use the result of Lemma~\ref{l:life-time-upper-pm}.

Observe that during the initial part of the life-time of a
protocluster the leading particle moves at rate one exchanging
positions with vacancies until at time $t_1$ it starts interacting
with the root. After this whence another particle from the body of
the protocluster reaches the leading one at time $t_2>t_1$ the
cluster is formed and one can apply Lemma~\ref{l:life-time-upper-pm}
to get the result. Therefore we are concerned only with these two
initial parts of the dynamics of a protocluster mentioned above.

During the time $0\le t< t_1$ the leading particle is moving
freely exchanging its position with vacancies.
Therefore %
$$ W(x^{t+1}[j(t),m(t+1)])=W(x^{t}[j(t),m(t)])-1 ,$$
because a new vacancy is taken into account while all already
present particles and vacancies contribute the same values as
before.

Observe that the site $j(t)$ cannot be occupied by a positive
particle, otherwise $$W(x^t[j(t),m(t)])\ge\tau+1$$ and hence
$$W(x^t[j(t)-1,m(t)])>0$$ which contradicts to the definition of
$j(t)$. Thus the boundary $j(t)$ moves at least at rate one. On the
other hand, the particle $\xi$ also cannot move faster than at rate
one. Therefore $\xi(0)\le j(0)$ implies that $\xi(t)\le j(t)$ during
the initial part of the life-time of a protocluster.

During the time $t_1\le t< t_2$ the leading particle stays at the
same site and thus there are no new vacancies. However, the state
of the leading particle is increasing at rate one which plays
exactly the same role and leads to the same result as in the
previous case. \qed

To deal with the jammed phase we need to consider the lower bound
for the BA. To this end we apply the functional $W'$ instead of
$W$ to show that the corresponding boundary is located inside of
the BA and that all particles from the corresponding segment must
join the cluster during its life-time.

\begin{lemma}\label{l:life-time-cluster-pm}
BA$(J^t)\supseteq x[j'(t),m(t)]$.
\end{lemma}
\proof The argument here is very similar to the one used in the
proof of Lemma~\ref{l:life-time-proto-pm} with the only difference
that we need to show now that for any particle $\xi(t)$ initially
located in the segment $x[j(0),m(0)]$ we have $\xi(t)\ge j(t)$. In
other words, that the boundary $j(t)$ cannot outran the particle
$\xi(t)$.

The state of the current leading particle of the cluster is
increasing by one until it will exchange positions with the root
particle. Therefore the value of the functional $W'$ is decreasing
at rate one. Observe now that the site $j'(t)$ cannot be occupied by
a positive particle (otherwise $W(x^t[j'(t),m(t)])\ge\tau)$ and
hence $W'(x^t[j'(t)-1,m(t)])\ge0$ which contradicts to the
definition of $j'(t)$). Thus the boundary $j'(t)$ moves at least by
rate one. On the other hand, the particle $\xi$ also cannot move
faster than at rate one. Therefore $j(t)$ being initially smaller
than $\xi(t)$ may outran the latter only if at some time $t'$ it
makes a jump longer than one. This can happen only if there is $i$
such that $j(t'-1)<i<m(t'-1)$ and $W'(x^t[i,m(t'-1)])=1$. In this
case one might expect that on the next time step we have
$W'(x^{t'}[i,m(t')])=0$. However this happens only if the segment
$x^{t'}[i-1,i]$ at time $t'-1$ is occupied only by vacancies, which
in turn means that $W'(x^{t'}[i-1,m(t'-1)])=0$ and hence
$i-1=j(t')$. Thus long boundary jumps cannot take place. \qed

This finishes the proof of Theorem~\ref{t:BA-pm}.

\bigskip

The functional $W'$ is everywhere positive in the case of an
infinite BA and is used to calculate the  second critical density
corresponding to the beginning of the jammed phase:
$$ \rho\tau - \t\rho\tau - (1-\rho-\t\rho) > 0
   \quad \Longrightarrow \quad
   \rho(1+\tau) - \t\rho(\tau-1) > 1
   \quad \Longrightarrow \quad
   \rho'_c=\frac{1+\t\rho(\tau-1)}{1+\tau} .$$
In the ``intermediate phase'' (see Section~\ref{s:hysteresis}) and
during the transient period of the free phase there are no infinite
life-time clusters but instead there might be {\em recurrent} ones,
which are appearing and disappearing again (in the intermediate
phase) and finite life-time clusters even in the eventually free
phase. To calculate the corresponding critical density $\rho_c$ we
need take into account that the duration of the long interaction is
$\tau+1$ instead of $\tau$. Therefore it can be calculated using the
functional $W$ in which the weight of a positive particle is
$(\tau+1)$ but the weight of a negative particle is $-\tau$. The
condition $W(x[m-n,m])\le0$ for some arbitrary large $n$ gives
$$ \rho(\tau+1) - \t\rho\tau - (1-\rho-\t\rho) > 0
   \quad \Longrightarrow \quad
   \rho(2+\tau) - \t\rho(\tau-1) > 1
   \quad \Longrightarrow \quad
   \rho_c=\frac{1+\t\rho(\tau-1)}{2+\tau} .$$
%Compare this relation to the 2nd type of freely moving hysteresis
%configuration.

\section{Duality relations and calculation of average velocities}
\label{s:duality}

In this Section we discuss connections between dynamics of
particles of different types.

A configuration $(-x)\in X$ is said to be {\em dual} to $x$. This is
indeed the case since the state of each particle is inverted from
positive to negative and vice versa. It is straightforward to check
that the dynamics of configurations $x^t$ and $-x^t$ are exactly the
same except for the direction of the particles movement. This type
of relations we call {\em duality} and they allow to simplify
significantly the analysis of such systems. In particular we shall
mainly discuss properties of positive particles having in mind that
the corresponding properties of negative particles can be obtained
by duality relations.

For a configuration $x\in X$ with densities of positive and negative
particles $\rho(x)$ and $\t\rho(x)$ respectively denote by
$V(\rho,\t\rho)$ the average velocity (if it exists) of positive
particles and by $\t{V}(\rho,\t\rho)$ the modulus of the average
velocity of negative particles. By the duality between positive and
negative particles we have $V(\rho,\t\rho)=\t{V}(\t\rho,\rho)$.

It turns out that the assumption about the existence of average
velocities together with some information about the structure of
configurations allows to calculate the average velocities in terms
of particle densities.

\begin{lemma}\label{l:duality}
Assume that in a configuration $x$ the densities and average
velocities of both types of particles are well defined and that
positive particles move freely (i.e. except from
the interactions with negative particles). Then %
\beq{e:duality}{V(\rho,\t\rho)
 = 1- (\tau-1)\t\rho~(V(\rho,\t\rho) + \t{V}(\rho,\t\rho))}%
\end{lemma}
\proof During the time $t>0$ a positive particle will meet
$(V(\rho,\t\rho)+\t{V}(\rho,\t\rho))\cdot t\t\rho + o(t)$ negative
ones. The notation $o(t)$ stands for a term which grows slower
than $t$, i.e. $o(t)/t\toas{t\to\infty}0$. Thus taking into
account that according to the assumption positive particles never
wait in clusters (which otherwise will give an additional
contribution) we have:
$$ V(\rho,\t\rho)\cdot t
 = t - (\tau-1)(V(\rho,\t\rho)
   + \t{V}(\rho,\t\rho))\cdot t \t\rho + o(t).$$
Dividing both sides of this relation by $t$ and passing to the
limit as $t\to\infty$ we get the result. \qed

\begin{lemma}\label{l:vel-free}
Assume additionally to the assumptions of Lemma~\ref{l:duality}
that negative particles are free also. Then
\beq{e:vel-free}{V(\rho,\t\rho)
 =\frac{1+(\tau-1)(\rho-\t\rho)}
                   {1+(\tau-1)(\rho+\t\rho)} .}
\end{lemma}
\proof By Lemma~\ref{l:duality} we have %
\bea{ V(\rho,\t\rho)
 \a= 1- (\tau-1)\t\rho~(V(\rho,\t\rho) + \t{V}(\rho,\t\rho)) \cr
   \t{V}(\rho,\t\rho)
 \a= 1- (\tau-1)\rho~(V(\rho,\t\rho) + \t{V}(\rho,\t\rho)) .}%
Denote $v= V(\rho,\t\rho) + \t{V}(\rho,\t\rho)$. Then adding the
above relations we get $v=2-(\tau-1)(\rho+\t\rho)v$, which yields
$$ v = \frac2{1+(\tau-1)(\rho+\t\rho)} $$
and hence
$$ V(\rho,\t\rho)
 = \frac{1+(\tau-1)(\rho - \t\rho)}{1+(\tau-1)(\rho + \t\rho)} .$$
In particular,
$$ V(\rho,\rho) = \frac1{1+2(\tau-1)\rho} .$$
\qed

\begin{lemma}\label{l:vel-jammed} Assume now that all negative
particles are free but positive particles form infinite life-time
clusters. Then $V(\rho) = \frac1\tau~(\frac1\rho - 1)$.
\end{lemma}
\proof According to our assumptions after at most $\tau+1$ time
steps the leading particles in the infinite life-time cluster will
take part only in short interactions with the root of the cluster.
Consider now the flux of positive particles through the root of the
cluster. According to the argument above after at most $\tau+1$
iterations the flux become equal to $1/\tau$ (i.e. after each $\tau$
moments of time a leading particle from the cluster body exchanges
positions with the cluster root). Thus $1/\tau = \rho(V(\rho,\t\rho)
+ \t{V}(\rho,\t\rho))$ and hence $V(\rho) = \frac1\tau~(\frac1\rho -
1)$. \qed

Observe that the difference between short and long interactions does
not matter both in the free and jammed phases because in the former
it does not change the timing, while in the latter in the steady
state it cannot happen if the length of the jam is greater than 1
(there are no gaps between particles in the cluster body).

\smallskip

The derivation of the average velocity in the free phase is based on
the assumption that $\t{V}(\rho,\t\rho)$ describes the flow of
noninteracting particles. If the density of negative particles is
becoming large enough then they start interact between themselves.
Assume that this density is so large that the corresponding velocity
belongs to the jammed
region and hence %
\bea{ V(\rho,\t\rho) \a= 1 - (\tau-1)(V(\rho,\t\rho)
                          + \t{V}(\rho,\t\rho))\t\rho \cr
  \t{V}(\rho,\t\rho) \a= \frac1\tau~(\frac1{\t\rho} -1) .}%
Then $V(\rho,\t\rho) = 1/\tau$ which does not depend on the particle
density and corresponds to the case when all the time a positive
particle interacts with negative ones, i.e. the negative particles
form an infinite life-time cluster. Clearly this means that if a
configuration belongs to $\Free_+^\infty\cap\Jam_-^\infty$ then the
average velocity of positive particles is equal to $1/\tau$.

\section{Proof of Theorems~\ref{t:main-pm} and \ref{t:phase-tr}}
\label{s:proofs}

These results will be proven more or less simultaneously through a
series of technical Lemmas.

\begin{lemma}\label{l:low-den-pm} Let a configuration
$x$ satisfies the condition %
\beq{e:low-den-pm}{(\tau+2)\rho(x)-(\tau-1)\t\rho(x)<1.}%
Then only finite life-time (proto)clusters may be present in the
configuration $x$. Moreover, there exists a partition of the
integer lattice into nonoverlapping BAs of these finite life-time
(proto)clusters and their complements.
\end{lemma}
\proof Assume on the contrary that there is a (proto)cluster
$J=x[m',m]$ with an infinite BA. By Theorem~\ref{t:BA-pm} and the
definition of the functional $W$ for any integer $n>1$ we have
$W(x[m-n,m])>0$. Rewriting this inequality in terms of particle
densities we get %
\bea{ q(x_m) \a+ (\tau+1)n\rho(x,[m-n,m-1])
          - \tau n \t\rho(x,[m-n,m-1]) \\
          \a- n(1 - \rho(x,[m-n,m-1]) - \t\rho(x,[m-n,m-1]))
          > 0 .} %
Dividing by $n$ and passing to the limit as $n\to\infty$ we obtain
$$ (\tau+2)\rho(x)-(\tau-1)\t\rho(x)\ge1 ,$$
which contradicts to the assumption (\ref{e:low-den-pm}).

To prove the second claim we need to demonstrate that the finite
BAs cannot be infinitely nested. The latter means that there
exists a sequence of (proto)clusters $\{J_n\}_n$ and an integer
$N$ such that for any $n>N$ we have BA$(J_n)\supset$BA$(J_N)$. If
this would be the case a positive particle located initially to
the right of $J_N$ will never become free despite the amount of
time it spends in each of the clusters is finite.

Assume that such an infinitely nested sequence of (proto)clusters
$\{J_n\}_n$ does exist. Then for some $m\in\IZ$ and any integer
$n>1$ by Theorem~\ref{t:BA-pm} we have $W(x[m,m+n])>0$. Rewriting
this in terms of particle densities we get %
\bea{ q(x_{m+n}) \a+ (\tau+1)n\rho(x,[m,m+n-1])
          - \tau n \t\rho(x,[m,m+n-1]) \\
          \a- n(1 - \rho(x,[m,m+n-1]) - \t\rho(x,[m,m+n-1]))
          > 0 .} %
Exactly as in the previous case we divide this inequality by $n$
and after passing to the limit as $n\to\infty$ we come to the
contradiction to the assumption (\ref{e:low-den-pm}). \qed

\n{\bf Remark}. The importance of the 2nd claim of
Lemma~\ref{l:low-den-pm} is that it might be possible that all
clusters have only finite life-times but a particle is joining
them one after another and is never becoming free.

\begin{lemma}\label{l:phase1-V-pm} Let a configuration
$x$ satisfy the conditions of Lemma~\ref{l:low-den-pm}. Then any
positive particle after some initial transient period will stop
interacting with other particles of the same sign (i.e. will
become free).
\end{lemma}
\proof By the previous Lemma there exists a partition of the
integer lattice into nonoverlapping BAs of finite life-time
(proto)clusters and their complements. Choose one  of these BAs
and consider a particle $\xi$ located originally immediately
before this BA. By the definition of the BA the particle $\xi$
does not belong to any BA of (proto)clusters located ahead of it
and hence it will never join any cluster of positive particles.
\qed

Results of Lemmas~\ref{l:low-den-pm}, \ref{l:phase1-V-pm}
prove claims (a,b) of Theorem~\ref{t:phase-tr}.

\begin{lemma}\label{l:high-den-pm} Let a configuration
$x$ satisfy the condition %
\beq{e:high-den-pm}{(\tau+1)\rho(x)-(\tau-1)\t\rho(x)>1.}%
Then for each $N\in\IZ$ there is an infinite life-time
(proto)cluster in $x$ located to the right from the site $N$.
\end{lemma}
\proof Fix some $N\in\IZ$ and assume on the contrary that any
positive particle in $x$ located at a site $i>N$ is either free or
belongs to a finite life-time (proto)cluster.

Thus there exists a sequence of finite BAs of (proto)clusters
$J_k$ such that BA$(J_k)\supseteq[j'_k,m_k]$ and
$j'_k<m_k<j'_{k+1}<m_{k+1}\toas{k\to\infty}\infty$. By the
definition of $j'_k$ for any $\ell>1$ we have
$$ \sum_{k=1}^\ell W'(x[j'_k,m_k]) < 0 .$$
Between the segments $[j'_k,m_k]$ there might be some additional
positive particles not belonging to any BA and thus each of them
satisfies the condition that either the number of vacancies
immediately ahead of it exceeds $\tau$ or it is preceded by a
negative particle. In both cases (as well if there are additional
negative particles) the value of the functional $W'$ on the
complete segment $[j'_1,m_\ell]$ instead of the union of segments
$[j'_k,m_k]$ is again negative.

Denoting $m:=j'_1,~n_\ell:=m_\ell-j'_1$ and rewriting the
inequality $W'(x[m,m+n_\ell-1])\le0$ in terms of particle
densities we get %
\bea{ q(x_{m+n_\ell}) \a+ (\tau+1)n\rho(x,[m,m+n_\ell-1])
          - \tau n \t\rho(x,[m,m+n_\ell-1]) \\
          \a - n(1 - \rho(x,[m,m+n_\ell-1])
                - \t\rho(x,[m,m+n_\ell-1])) < 0 .} %
Dividing this inequality by $n$ and passing to the limit as
$\ell\to\infty$ we obtain
$$ (\tau+1)\rho(x)-(\tau-1)\t\rho(x)\le1 .$$
Thus we came to the contradiction to the assumption
(\ref{e:high-den-pm}). \qed

Claims (c,d) of Theorem~\ref{t:phase-tr} follow from this result.

\bigskip

In the hysteresis phase the protoclusters {\em might} be present,
which makes the main difference to the free phase. Moreover, it
might be possible that despite the protoclusters have only finite
life-times their BAs are infinitely nested which prevents positive
particles to become free eventually.

\bigskip

Now we are ready to prove Theorem~\ref{t:main-pm}.

We start with the case $(\rho(x),\t\rho(x))\in A$. Then the
particle densities satisfy the system of inequalities %
\bea{\a(\tau+2)\rho<1+(\tau-1)\t\rho ,\cr
     \a(\tau+2)\t\rho<1+(\tau-1)\rho .}%
By Lemma~\ref{l:phase1-V-pm} and the first of these inequalities all
positive particles will eventually become free. On the other hand,
applying the duality relation together with the second inequality,
by Lemma~\ref{l:phase1-V-pm} we get the same property for negative
particles. Therefore exact relations for the average velocities
$(V(x),\t{V}(x))$ follow from Lemma~\ref{l:vel-free}.

Assume now that $(\rho(x),\t\rho(x))\in C$. Then %
\bea{\a(\tau+1)\rho>1+(\tau-1)\t\rho ,\cr
     \a(\tau+2)\t\rho<1+(\tau-1)\rho .}%
The first of these inequalities implies by Lemma~\ref{l:high-den-pm}
that in the configuration $x$ there are infinitely many infinite
life-time (proto)clusters of positive particles, while from the
second inequality it follows that all negative particles will
eventually become free. Thus some negative particle will become a
short root of an infinite life-time cluster and from that moment of
time will move exactly by one site in $\tau$ time steps (exchanging
positions with positive particles from the body of the cluster). Now
since the average velocity is well defined for a single particle, by
Lemma~\ref{l:velocity-inv-pm} the same result holds for all
particles of the same sign and thus $\t{V}(x)=1/\tau$. Once we get
this we obtain the relation for the average velocity for positive
particles by the duality relation.

The case $(\rho(x),\t\rho(x))\in B$ is considered similarly to the
previous one, except for in this case positive particles get the
average velocity $1/\tau$.

Theorems~\ref{t:main-pm} and \ref{t:phase-tr} are proven. \qed

\section{Hysteresis region}\label{s:hysteresis}

Our previous results are concerned with the behavior of
configurations in ``pure'' free and jammed phases. Now we shall show
that in the region $H$ located between these phases (see
Fig.\ref{f:fund-diag-pm}(right)) the coexistence of free and jammed
configurations takes place. In other words in this region there
might be both freely moving and jammed configurations with the same
particle densities. Due to the duality between positive and negative
particles in what follows we consider only the case when
$\rho\ge\t\rho$.

\begin{lemma}\label{l:hyst-free} Configurations from the set
$\Free_+^\infty$ are dense in %
$H_+:=\{(\rho,\t\rho)\in H:~~(\tau+1)\rho\le1+(\tau-1)\t\rho\}$.
\end{lemma}

\proof Consider a family of spatially periodic configurations
parameterized by nonnegative integers $n,m,k_1,\dots,k_n$ (only
one spatial period is shown): %
\beq{e:free1}{\overbrace{
  +\underbrace{\cdots}_{\tau+2k_1}
  +\underbrace{\cdots}_{\tau+2k_2} \dots
  +\underbrace{\cdots}_{\tau+2k_n} }^{(\tau+1)n+2N}
  \overbrace{
  \underbrace{+-}{}\underbrace{+-}{}\dots\underbrace{+-}{}
  }^{2m} } %
Setting $N:=\sum_{i=1}^nk_i$ the length of the period can be written
as $(\tau+1)n+2(m+N)$. The value $\tau+2k_i$ is equal to the
duration of the short interaction plus an arbitrary even number.
This guarantees that only short interactions may take place under
dynamics. On the other hand, the number of vacancies between
consecutive positive particles is not smaller than $\tau$, which
implies that for any $k_i$ a positive particle cannot catch up with
another positive particle (since only short interactions may take
place). Therefore $x\in\Free_+$. %
%Similarly one shows that $x\in\Free_-$.

Let a configuration $x$ belongs to the family~(\ref{e:free1}). Then %
\bea{
   \rho(x)\a=\frac{n+m}{(\tau+1)n+2(m+N)} =:r(n,m,N) ,\\
   \t\rho(x)\a=\frac{m}{(\tau+1)n+2(m+N)} =:\t{r}(n,m,N) .} %
Observe now that for any $u,v,z\ge0$ such that $u+v>0$ %
the identities %
$$ r(u,v,z)\equiv r(\frac{u}{u+v},\frac{v}{u+v},\frac{z}{u+v}),
\qquad
   r(u,v,z)\equiv r(\frac{u}{u+v},\frac{v}{u+v},\frac{z}{u+v}) $$%
hold true and that %
$$ 1+(\tau-1)\t{r}(u,v,z) - (\tau+1)r(u,v,z)
 = \frac{2z}{(\tau+1)u+2(v+z)} \ge 0 $$ %
Therefore the ``right'' boundary of the region $H$ (corresponding to
the case $z=0$) contains configurations from $\Free_+$. Moreover,
the points on the ``right'' boundary for which the above
representation takes place are dense. Additionally we see that these
functions decay with respect to each their argument. Therefore for a
given pair or real numbers $u,v\ge0,~u+v>0$ the pair of functions
$(r(u,v,z),~\t{r}(u,v,z))$, considered as functions of the third
real argument $z$, parameterizes a curve
$R_{u,v}(z):=(r(u,v,z),~\t{r}(u,v,z))$ which starts (when $z=0$) at
some point on the ``right'' boundary of $H$ and goes through $H$
until it reaches the ``left'' boundary.

For $d\ge1$ denote
$\dist_d((u_1,\dots,u_d),(u_1',\dots,u_d')):=\sum_i|u_i-u_i'|$. Then %
\beq{e:dist_curve}{\dist_2(R_{u,v}(z),R_{u',v'}(z'))
  \le \frac{2\tau\dist_3((u,v,z),(u',v',z'))}
           {(u+v+z)(u'+v'+z')} } %
for any pair of triples $(u,v,z),~(u',v',z')$. Therefore the
curves $R_{u,v}$ with $u,v\in\IR_+\cup\{0\}$ fill in the entire
region $H_+$.

It remains to show that for any real $(\rho,\t\rho)\in H_+$ and
any $\ep>0$ there exists a triple of nonnegative integers
$(n,m,N)$ such that %
$$ \dist_2((\rho,\t\rho),~R_{n,m}(N)) < \ep .$$ To do this set
$n':=\frac{K\rho}{\rho+\t\rho}, ~m':=\frac{K\t\rho}{\rho+\t\rho}$.
Then for any $K>0$ the curve $R_{n',m'}$ goes through the point
$(\rho,\t\rho)$ at some real parameter value $N'$. Since points
on the ``right'' boundary of $H$ represented by integer values of
$(n,m,N)$ are dense we may choose two sequences of nonnegative
integers $\{n_k\}, \{m_k\}\toas{k\to\infty}\infty$ such that $$
\dist_2(R_{n',m'}(0),~R_{n_k,m_k}(0))\toas{k\to\infty} 0 .$$

Therefore setting %
$$ u':=\frac{\rho}{\rho+\t\rho}, ~v':=\frac{K\t\rho}{\rho+\t\rho},
~ z':=\frac{N'}{n'+m'}; $$ %
$$ u_k:=\frac{n_k}{n_k+m_k}, ~v_k:=\frac{m_k}{n_k+m_k},
~ z_k:=\frac{N}{n_k+m_k}; $$ %
we get
$$ \dist_2((u_k,v_k),(u',v')) \toas{k\to\infty} 0 .$$ Using
(\ref{e:dist_curve}) we obtain %
\bea{ \dist_2(R_{u_k,v_k}(z_k),R_{u',v'}(z'))
   \a\le \frac{2\tau\dist_3((u_k,v_k,z_k),(u',v',z')}{(1+z_k)(1+z')}
   \\
   \a\le 2\tau(\dist_2((u_k,v_k),(u',v')) + |z_k - z'| ) .} %
Therefore the curves $R_{n_k,m_k}$ come arbitrary close to the point
$(\rho,\t\rho)$ as $k\to\infty$ at some integer
parameter values $N$, for which %
$$|z_k-z'|\toas{k\to\infty} 0$$ %
and this finishes the proof. \qed

To study the ``left'' boundary of $H$ consider another family
of spatially periodic configurations: %
\beq{e:free2}{\overbrace{
  +\underbrace{\cdots}_{\tau+1+2k_1}
  +\underbrace{\cdots}_{\tau+1+2k_2} \dots
  +\underbrace{\cdots}_{\tau+1+2k_n} }^{(\tau+2)n+2N}
  \overbrace{
  \underbrace{+\cdot-}{}\underbrace{+\cdot-}{}
  \dots\underbrace{+\cdot-}{}
  }^{3m} } %
The notation here is the same as for the family~(\ref{e:free1})
and the length of the period is equal to $(\tau+2)n+2(m+N)$. A
direct check shows that again in this family both positive and
negative particles are free. The difference is the absence of long
interactions under dynamics for configurations from
(\ref{e:free1}) and their inevitable presence for configurations
from (\ref{e:free2}) for an arbitrary large time. If a
configuration $x$ belongs to the family~(\ref{e:free2}) then
$$ \rho(x)=\frac{n+m}{(\tau+2)n+2(m+N)} \qquad
   \t\rho(x)=\frac{m}{(\tau+2)n+2(m+N)} .$$
Setting again $N=0$ we obtain the ``left'' boundary of $H$, namely
$$ (\tau+2)\rho(x) \equiv 1+(\tau-1)\t\rho(x) .$$

On the other hand, for the first family
$$1\ge(\rho-\t\rho)(\tau+1)+2\t\rho
  \quad \Longrightarrow \quad
  \rho\le\frac{1+\t\rho(\tau-1)}{\tau+1}$$
while for the second family
$$1\ge(\rho-\t\rho)(\tau+2)+3\t\rho
  \quad \Longrightarrow \quad
  \rho\le\frac{1+\t\rho(\tau-1)}{\tau+2} .$$
Comparing these inequalities with the critical densities $\rho'_c$
and $\rho_c$ we see that the family (\ref{e:free2}) belongs to the
free phase while the family~(\ref{e:free1}) does not belong there.

\bigskip

It remains to show that in this region there are jammed (for all
times) configurations as well. According to
Lemma~\ref{l:low-den-pm} the inequality~(\ref{l:low-den-pm})
(defining the ``right'' boundary of $H_+$) guarantees the absence
of infinite life-time (proto)clusters in $H$. This, together with
the denseness of densities corresponding to free configurations in
$H$ (see Lemma~\ref{l:hyst-free}) explains the subtlety of this
problem: one needs to find configurations having only
(proto)clusters of arbitrary large but finite life-time. To
understand how this can happen consider the dynamics of a
spatially periodic configuration $x\in X$ with $\tau=2$:

\?{ %commented old version
(numbers 2 and 3 indicate
positions of positive particles with states greater than one): %
\def\cd{\cdot} %
\bea{
   \a +  \cd\cd\cd  +  -  +  - \\ %
   \a \cd  +  \cd\cd  2  -  2  - \\ %
   \a \cd\cd  +  \cd  -  +  -  + \\ %
   \a   +  \cd  2  \cd  -  2  -  \cd \\ %
   \a \cdot  +    3  \cdot  -  -  +  \cdot \\ %
   \a \cdot  +    -  \cdot  +  -  \cdot +  \\ %
   \a   +    2    -  \cdot  2   - \cdot\cdot \\ %
   \a   +    -    +  \cdot  -   + \cdot\cdot \\ %
   \a   2    -    2  \cdot  - \cdot  + \cdot \\ %
   \a   -    +    3  \cdot  - \cdot\cdot + \\ %
   \a   -    +    -  \cdot  + \cdot\cdot 2 \\ %
   \a   +    2    -  \cdot\cdot +  \cdot - }%
}

\def\vb{\vskip-0.65cm}
$$ 10001\t11\t1 $$\vb
$$ 01002\t22\t2 $$\vb
$$ 0010\t11\t11 $$\vb
$$ 1020\t22\t20 $$\vb
$$ 0130\t3\t110 $$\vb
$$ 01\t101\t101 $$\vb
$$ 12\t202\t200 $$\vb
$$ 1\t110\t1100 $$\vb
$$ 2\t220\t2010 $$\vb
$$ \t1130\t3001 $$\vb
$$ \t21\t101002 $$\vb
$$ 12\t20010\t1 $$

\bigskip

Using the notation introduced in the proof of
Lemma~\ref{l:hyst-free} we have $n=3,m=2,N=3$. Thus
$r(n,m,N)=3/8$, $\t{r}(n,m,N)=1/4$ and hence %
$$ (\tau+1)r(n,m,N) = 9/8 < 1 + (\tau-1)\t{r}(n,m,N) = 5/4
 < (\tau+2)r(n,m,N) = 3/2 ,$$
i.e. $x$ belongs to the region $H$. During the dynamics each
positive particle in turn\footnote{The sequence
   of types of interactions for a single particle is periodic but
   may be arbitrary.} %
makes short and long interactions with negative ones. Finite
(proto)clusters are present for all time. On the other hand only
trivial life-time clusters are present (at moments of time when
one of positive particles catches up with another), hence the
configuration in this example does not belong to the jammed phase.
Nevertheless at any moment of time there are infinitely nested BAs
of trivial protoclusters.

Our numerical simulations, which were carried out for a broad
region of spatial periods $L$ (up to $L=200$) indicate that for
each spatial period $L$ for all pairs of positive integers $n,m$
for which $n+m<L$ and $(r(n,m,L-n-m),\t{r}(n,m,L-n-m))\in H$ there
were configurations with only trivial (proto)clusters for
arbitrary large time. This allows to formulate the following

\bigskip

\n{\bf Hypothesis}. {\em For each pair of densities belonging to the
region $H$ there exists an admissible configuration having positive
and negative particles with these densities for which only trivial
(proto)clusters are present for arbitrary large time.}

\bigskip

For spatially periodic configurations the existence of average
velocities follows from the obvious time-periodicity of the dynamics
of particles. However in the general case this is not proven.
Moreover, without the assumption about the existence of particle
densities one can construct configurations for which average
velocities do not exist.

\section{Dynamics of an active tracer} \label{s:active}

To model the motion of an active tracer in the traffic flow we add
to the flow a single particle (representing the tracer). If the
tracer moves in the same direction as the flow we think about this
particle as having the same sign as particles of the flow, and
having the opposite sign when the tracer moves against the flow. In
distinction to the passive tracer studied in \cite{Bl3} the active
one changes the behavior of the flow due to interactions with other
particles. Clearly a single positive particle added to the flow of
positive particles do not change its asymptotic behavior. On the
other hand the case when the tracer moves against the flow is
equivalent to the model discussed in this paper with zero density of
positive particles (tracer) and a certain density of negative
particles (the flow). According to Theorem~\ref{t:main-pm} we get
the following average velocity of the active tracer as a function of
the density of the flow $\rho$:

\bea{%e:active-tracer-pm}{
V_+(\rho) \a= \function{1       &\mbox{if } 0\le\rho\le1/2 \\
                      1/\rho  &\mbox{otherwise }} \cr
V_-(\rho) \a= \function{\frac{1-(\tau-1)\rho}{1+(\tau-1)\rho}
                              &\mbox{if } 0\le\rho<\frac1{\tau+2} \\
                         \frac1\tau
                              &\mbox{if } \frac1{\tau+1}<\rho\le1 } }%
Here the average velocity $V_+(\rho)$ corresponds to the case when
the tracer moves in the same direction as the flow, while
$V_-(\rho)$ describes the case when the tracer moves in the
opposite direction. Fig.~\ref{f:active-tracer-pm} shows that there
is a region of flow densities when it is more advantageous to move
against the flow than along it. The explanation is that in high
density region the time necessary to exchange positions with a
particle moving in the opposite direction is becoming smaller than
the time a particle from the flow spends in inevitably being
present large clusters.

%%%%%%%%%%%%%%%%%%%%%%%%%%%%%%%%%%%%%%%%%%%%%%%%%%%%%%%%%%%%%%%%%%
%% Dynamics of an active tracer (+-)
\Bfig(150,120)
      {\footnotesize{
       %% picture (a)
       \put(0,0){\vector(1,0){150}} \put(0,0){\vector(0,1){120}}
       \thicklines
       \bezier{100}(0,105)(35,105)(70,105)%phase 1
       \bezier{200}(70,105)(100,25)(135,0) %phase 3
       %\bezier{10}(46,80)(54,86)(58,75)  %phase 2
       %\bezier{10}(58,75)(64,55)(70,70)  %---"---
       \thinlines
       %\bezier{10}(46,80)(55,85)(70,90) %continuation of phase 1
       %\bezier{10}(63,87.5)(67,80)(70,70) %continuation of phase 3
       \bezier{25}(30,110)(30,55)(30,0) %1st boundary (46,78)
       \bezier{25}(40,110)(40,55)(40,0) %2nd boundary
       \bezier{25}(70,110)(70,55)(70,0) %1/2
       \bezier{100}(0,105)(15,100)(30,90) %V_ (-)-particles
       \bezier{200}(40,30)(90,30)(135,30) %V_ (-)-particles
       \bezier{10}(0,30)(20,30)(40,30) %V_ (-)-particles
       \put(15,55){$V_-$} \put(90,80){$V_+$} %(-)-particles
       \put(18,-10){$\frac1{\tau+2}$} \put(37,-10){$\frac1{\tau+1}$}
       \put(145,-8){$\rho$} \put(130,-8){1} \put(67,-10){$\frac12$}
       \put(-10,115){$V$}   \put(-10,105){$1$}
       \put(-10,27){$\frac1\tau$}
      }}
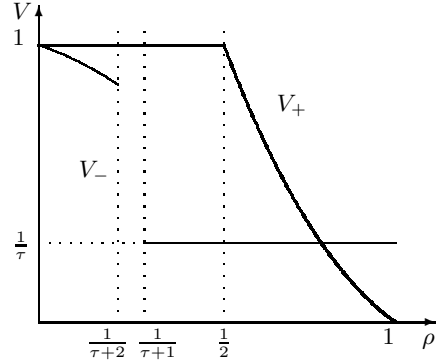
{Dependence of the average velocity of the active tracer $\V_\pm$
on the particle density $\rho$. $\V_+$ -- along the flow (thick
line), $\V_-$ -- against it (thin line).
\label{f:active-tracer-pm}}
%%%%%%%%%%%%%%%%%%%%%%%%%%%%%%%%%%%%%%%%%%%%%%%%%%%%%%%%%%%%%%%%%

\section{Configurations having no particle densities}
\label{s:generalization}

So far we were discussing statistical properties of configurations
for which both positive and negative particle densities are well
defined. In general one or both of these densities might be not well
defined and one needs to deal with lower and upper densities. Recall
that Lemma~\ref{l:density-preservation} claims their invariance with
respect to dynamics. The following result extends the description of
the Phase diagram for this more general case.

\begin{theorem}\label{t:phase-tr-pm} Let $x\in X$. If \\%
(a) $(\tau+2)\rho^+(x)<1+(\tau-1)\t\rho^-(x)$ then
    $x\in\Free_+^\infty$.\\
(b) $(\tau+2)\t\rho^+(x)<1+(\tau-1)\rho^-(x)$ then
    $x\in\Free_-^\infty$.\\
(c) $(\tau+1)\rho^-(x)>1+(\tau-1)\t\rho^+(x)$ then
    $x\in\Jam_+^\infty$.\\
(d) $(\tau+1)\t\rho^-(x)>1+(\tau-1)\rho^+(x)$ then
    $x\in\Jam_-^\infty$.
\end{theorem}

The proof is similar to the proof of Theorem~\ref{t:phase-tr} and
therefore we skip it.

In distinction to the case when the densities are well defined,
here we cannot prove the existence of average velocities. Instead
one considers upper and lower particle velocities for which it is
possible to get a representation of the same sort as in
Theorem~\ref{t:main-pm}.

\section*{Acknowledgement} The author would like to thank the referees
for valuable suggestions and comments.

%\footnotesize %%%%%%%%%%%%%%%%%%%%%%%%%%%%%%%%%%%%

\end{document}